\documentclass[a4paper,11pt]{article}
\usepackage{amssymb}
\usepackage{amsthm}
\usepackage{pstricks}
\usepackage{amsmath}
\usepackage{epsfig}
\usepackage{mathrsfs}

\def\XXint#1#2#3{{\setbox0=\hbox{$#1{#2#3}{\int}$}
     \vcenter{\hbox{$#2#3$}}\kern-.5\wd0}}

\theoremstyle{definition} \theoremstyle{remark}

\begin{document}

\title{\textbf{The perturbed Bessel equation, I. A Duality Theorem. }}
\author{\textbf{V.P. Gurarii and D.W.H. Gillam}}
\maketitle

\begin{abstract}
The Euler-Gauss linear transformation formula for the hypergeometric
function was extended by Goursat for the case of logarithmic singularities.
By replacing the perturbed Bessel differential equation by a monodromic
functional equation, and studying this equation separately from the
differential equation by an appropriate Laplace-Borel technique, we
associate with the latter equation another monodromic relation in the dual
complex plane. This enables us to prove a duality theorem and to extend
Goursat's formula to much larger classes of functions.
\end{abstract}

\section{Introduction}

We propose a new approach for the study of the perturbed Bessel differential
equation (pBde). This approach is based on the duality between the
functional monodromic relation generated by the pBde and its Laplace-Borel
dual. The solutions of the latter relation can be viewed as generalized
hypergeometric functions in the sense that they inherit a number of
important properties of the classical hypergeometric functions.

The pBde, with an infinite number of parameters, gives rise to a monodromic
relation with only two parameters
\begin{equation}
P\left( \zeta e^{\pi i}\right) =P\left( \zeta e^{-\pi i}\right) +Te^{-a\zeta
}P\left( \zeta \right) ,  \label{*}
\end{equation}
where $a>0,T$ is a complex constant which depends on all parameters of the
pBde, $P\left( \zeta \right) $ is analytic at every point $\zeta \neq 0$ and
bounded at infinity in any sectorial region. The dual function $F\left(
t\right) $ is analytic in the $t$-plane cut along the interval $\left(
-\infty ,-a\right) $ of the$\ $real line and admits an analytical
continuation along any path not containing points $t=0$ and $t=-a.$ Among
other properties which follow from (\ref{*}) are that the function $F\left(
t\right) $\ has exponential growth at infinity of minimal type in any
sectorial region, and satisfies the monodromic relation
\begin{equation}
F\left( \left( t+a\right) e^{\pi i}-a\right) =F\left( \left( t+a\right)
e^{-\pi i}-a\right) -TF\left( -\left( t+a\right) e^{\pi i}\right) .
\label{**}
\end{equation}

The essence\ of our approach is in considering both the monodromic relations
for the pBde and its dual quite separately from the differential equations
which generate them. In brief, given $a$ and $T$, we consider the linear
spaces $S$ and $H$ of all functions retaining the analytic properties stated
in the previous paragraph, and satisfying (\ref{*}) and (\ref{**}),
respectively, and prove the following result:

\textbf{Duality Theorem}. \textit{For fixed values of }$a$\textit{\ and }$T,$%
\textit{\ the Laplace transform operator }$\mathcal{L}:H\rightarrow S$%
\textit{\ is a bijection}$.$

The monodromic relation (\ref{**}) does not, on its own, determine the
behavior of $F\left( t\right) $ at $t=-a$. However, the relation (\ref{*}),
together with the conditions on $P\left( \zeta \right) $ stated above, shows
that $F\left( t\right) $ has a logarithmic singularity at $t=-a$. Moreover
we are able to prove the stronger statement:

\textbf{Theorem 2}. \textit{For fixed values of }$a$\textit{\ and }$T,\,$%
\textit{assume that }$F\left( t\right) $\textit{\ belongs to the linear
space H. Then in the region} $\left| t+a\right| <a,\left| \arg \left(
t+a\right) \right| <\pi ,$ \textit{the following relationship is valid} \
\begin{equation}
F\left( t\right) =-TF\left( -\left( t+a\right) \right) \log \left(
t+a\right) +\sum_{k=0}^\infty A_k\left( t+a\right) ^k,  \label{***}
\end{equation}
\textit{where the series in (\ref{***}) is absolutely convergent and }$A_k$%
\textit{\ are complex coefficients which can be evaluated explicitly.}

Formula (\ref{***}) can be viewed as an extension of the special case of the
Goursat linear transformation formula for the hypergeometric function $%
_2F_1\left( \mathfrak{a},\mathfrak{b};\mathfrak{c};t\right) $, valid when $%
\mathfrak{a}+\mathfrak{b}=\mathfrak{c}=1$\ and $_2F_1\left( \mathfrak{a},%
\mathfrak{b};\mathfrak{c};t\right) $\ has a logarithmic singularity at $t=1.$
While Goursat's approach relied on properties specific to the hypergeometric
function, we establish a linear transformation formula for all elements of
space $H.$

The above results allow us, for fixed values of $a$\ and $T,$ to provide all
functions from the linear space $S$ with asymptotic expansions in the $\zeta
$-plane, and using appropriate exponentially small terms to supply these
expansions with uniform error bounds, see \cite{GGM-2}. Returning to the
original differential equation and developing our technique further, we
derive an explicit formula for the connection coefficient $T$, see \cite%
{GGM-3}.

\section{Two Relations for the Hypergeometric Function}

In this section we consider the linear transformation formula and monodromic
relation for the classical hypergeometric function \textit{\ }$_2F_1\left(
\mathfrak{a},\mathfrak{b};\mathfrak{c};t\right) .$ The Euler linear
transformation formula which determines the behavior of the hypergeometric
function \textit{\ }$_2F_1\left( \mathfrak{a},\mathfrak{b};\mathfrak{c}%
;t\right) $ at point $t=1$ is well known (see for example \cite{Abramowitz},
\textbf{15.3.6}) and is valid for all values of the parameters except for
cases $\mathfrak{a}+\mathfrak{b}=\mathfrak{c}+m,$\ where $m$\ is an integer.
The associated formula
\begin{align}
& _2F_1\left( \mathfrak{a},\mathfrak{b};\mathfrak{c};1-\left( 1-t\right)
e^{\pm 2\pi i}\right) -\left( _2F_1\left( \mathfrak{a},\mathfrak{b};%
\mathfrak{c};t\right) \right)  \notag \\
=T^{\pm }(\mathfrak{a},\mathfrak{b},\mathfrak{c})& \left( 1-t\right) ^{%
\mathfrak{c}-\mathfrak{a}-\mathfrak{b}}\,_2F_1\left( \mathfrak{c}-\mathfrak{a%
},\mathfrak{c}-\mathfrak{b};\mathfrak{c}-\mathfrak{a}-\mathfrak{b}%
+1;1-t\right) ,  \label{****}
\end{align}
where
\begin{equation}
T^{\pm }(\mathfrak{a},\mathfrak{b},\mathfrak{c})=\mp 2{\pi }ie^{\pm \pi
i\left( \mathfrak{c}-\mathfrak{a}-\mathfrak{b}\right) }\frac{\Gamma \left(
\mathfrak{c}\right) }{\Gamma \left( \mathfrak{a}\right) \Gamma \left(
\mathfrak{b}\right) \Gamma \left( \mathfrak{c}-\mathfrak{a}-\mathfrak{b}%
+1\right) }  \label{T}
\end{equation}
is less known and can be derived from the Euler linear transformation
formula (see for example \cite{Kratzer}, section (4.2.4) formula (25). We
refer to this formula as the \textit{monodromic relation for }$_2F_1\left(
\mathfrak{a},\mathfrak{b};\mathfrak{c};t\right) .$ It is interesting to note
that formula (\ref{****}) is valid for all values of the parameters
including cases $\mathfrak{a}+\mathfrak{b}=\mathfrak{c}+m$ referred to
above. For the particular case $\mathfrak{a}+\mathfrak{b}=\mathfrak{c}$,
this formula takes the form
\begin{align}
_2& F_1\left( \mathfrak{a},\mathfrak{b};\mathfrak{a}+\mathfrak{b};1-\left(
1-t\right) e^{\pm 2\pi i}\right) -\text{{}}_2F_1\left( \mathfrak{a},%
\mathfrak{b};\mathfrak{a}+\mathfrak{b};t\right)  \notag \\
=& \mp 2{\pi }i\frac{\Gamma \left( \mathfrak{a}+\mathfrak{b}\right) }{\Gamma
\left( \mathfrak{a}\right) \Gamma \left( \mathfrak{b}\right) }{}_2F_1\left(
\mathfrak{a},\mathfrak{b};1;1-t\right) .  \label{ltf-B}
\end{align}

In 1881, E. Goursat, \cite{Goursat}, \cite{Goursat-1}, when studying the
logarithmic singularities of the hypergeometric function \textit{\ }$%
_{2}F_{1}\left( \mathfrak{a},\mathfrak{b};\mathfrak{c};t\right) $, derived
supplementary linear transformation formulas for the cases $\mathfrak{a}+%
\mathfrak{b}=\mathfrak{c}+m$\ where $m$\ is an integer, which are quite
different from Euler formula given in \cite{Abramowitz}, \textbf{15.3.6}. In
particular, he proved that if $\mathfrak{a}+\mathfrak{b}=\mathfrak{c}$ then
for $|\arg (1-t)|<\pi $ and $|1-t|<1$ the following formula is valid, (\cite%
{Abramowitz}, \textbf{15.3.10,})
\begin{align}
& 2F_{1}\left( \mathfrak{a},\mathfrak{b};\mathfrak{a}+\mathfrak{b};t\right)
\notag \\
={}& \frac{\Gamma \left( \mathfrak{a}+\mathfrak{b}\right) }{\Gamma \left(
\mathfrak{a}\right) \Gamma \left( \mathfrak{b}\right) }\left(
\sum_{n=0}^{\infty }\frac{\left( \mathfrak{a}\right) _{n}\left( \mathfrak{b}%
\right) _{n}}{\left( n!\right) ^{2}}d(n)\left( 1-t\right) ^{n}-\log \left(
1-t\right) {}_{2}F_{1}\left( \mathfrak{a},\mathfrak{b};1;1-t\right) \right) ,
\label{log-*}
\end{align}%
where
\begin{equation}
d(n)=2\psi \left( n+1\right) -\psi \left( \mathfrak{a}+n\right) -\psi \left(
\mathfrak{b}+n\right) ,\psi =\left( \log \Gamma \right) ^{\prime }.
\label{d-n-*}
\end{equation}%
Note that formula (\ref{log-*}) also yields the monodromic relation (\ref%
{ltf-B}). But the task of deriving (\ref{log-*}) using only the monodromic
relation (\ref{ltf-B}) seems difficult, if not insurmountable, even for the
particular case when $\mathfrak{a}+\mathfrak{b}=1.$ In this latter case
formula (\ref{ltf-B}) takes the simpler form
\begin{align}
_{2}F_{1}& \left( \mathfrak{a},1-\mathfrak{a};1;1-\left( 1-t\right) e^{\pm
2\pi i}\right) -{}_{2}F_{1}\left( \mathfrak{a},1-\mathfrak{a};1;t\right)
\notag \\
=& T\left( \mathfrak{a}\right) {}_{2}F_{1}\left( \mathfrak{a},1-\mathfrak{a}%
;1;1-t\right) ,  \label{monrel}
\end{align}%
where
\begin{equation}
T\left( \mathfrak{a}\right) =\mp 2i\sin \left( \pi \mathfrak{a}\right) ,
\label{Ta}
\end{equation}
and the same function occurs in each term of the relation. Even for this
special case formula (\ref{log-*}) remains as complicated as before.

Given $a>0$, replacing $t\rightarrow -\frac{t}{a}$ in (\ref{monrel})$,$ this
relation for the hypergeometric function ${}_{2}F_{1}\left( \mathfrak{a},1-%
\mathfrak{a};1;-\frac{t}{a}\right) $ can be rewritten in the form (\ref{**}%
), where $T$ is given by (\ref{Ta}). In the next section we show that the
dual equation (\ref{*}) for (\ref{**}) arises from the Bessel differential
equation. It turns out that the monodromic relation generated by the pBde is
indistinguishable from the analogous relation generated by the Bessel
differential equation. Therefore by identifying additional minimal
properties of \textit{\ }$_{2}F_{1}\left( \mathfrak{a},\mathfrak{b};%
\mathfrak{c};t\right) $\ that are required to derive (\ref{log-*}) from (\ref%
{ltf-B}) for the case $\mathfrak{a}+\mathfrak{b}=\mathfrak{c=}1$, or
Bessel's case$,$ we find similar conditions for the monodromic relation for
the pBde.

While our Theorem 2 appears to involve the derivation of (\ref{log-*}) from (%
\ref{monrel}), assuming $\mathfrak{a}+\mathfrak{b}=1,$ it in fact provides a
partial answer to a deeper question regarding the additional properties of
\textit{\ }$_{2}F_{1}\left( \mathfrak{a},\mathfrak{b};\mathfrak{c};t\right) $%
\ required to derive (\ref{log-*}) from (\ref{ltf-B}). To answer this in
full the results presented here for the pBde would need to be extended to a
more general case of the perturbed Whittaker differential equation (pWde).

\section{\textbf{The }Spaces $S_{a,T}$ and $H_{a,T}$}

Given $a>0$ and a complex number $T$, we define two function spaces $S_{a,T}$
and $H_{a,T}$ in terms of the monodromic relations (\ref{*}) and (\ref{**}).

\subsection{Monodromic Relations generated by the pBde}

We begin with a differential equation of the form
\begin{equation}
\frac{d^2u}{d\zeta ^2}=\left( \frac{a^2}4+\frac{a_0}{\zeta ^2}+\frac{a_1}{%
\zeta ^4}+\ldots \right) u,  \label{T-Airy}
\end{equation}
where $a>0$, $a_0,a_1,\ldots $ are complex numbers, and the series is
convergent for any complex $\zeta \neq 0$. Thus, if we set
\begin{equation}
A\left( \zeta \right) =a_0+\frac{a_1}{\zeta ^2}+\ldots ,  \label{a-coef}
\end{equation}
then $A\left( \zeta ^{-1}\right) $ is an entire even function. Equation (\ref%
{T-Airy}) can be considered as a perturbation of the standard Bessel
equation. Indeed, upon applying the transformation $\zeta =iz,u=z^{\frac
12}y,$ Bessel's equation
\begin{equation}
y^{\prime \prime }+\frac 1zy^{\prime }+\left( 1-\frac{\nu ^2}{z^2}\right) y=0
\label{Bessel eq}
\end{equation}
reduces\ to the following particular case of (\ref{T-Airy}) with $a=2$ and $%
A\left( \zeta \right) =\nu ^2-\frac 14$

\
\begin{equation}
\frac{d^{2}u}{d\zeta ^{2}}=\left( 1-\frac{\frac{1}{4}-\nu ^{2}}{\zeta ^{2}}%
\right) u.  \label{Bessel}
\end{equation}

Clearly, (\ref{T-Airy}) is invariant under the rotations $\zeta \rightarrow
\zeta e^{\pm \pi i}$ and also under the reflection $\zeta \rightarrow -\zeta
.$ Since $\zeta =0$ and $\zeta =\infty $ are the only singular points of (%
\ref{T-Airy}), any solution of this differential equation admits an
analytical continuation in the $\zeta $-plane along any path not containing
these points.

Let $u\left( \zeta \right) $ be a solution of (\ref{T-Airy}) decaying on the
positive ray. Representing this solution in the form
\begin{equation}
u\left( \zeta \right) =e^{-\frac{a}{2}\zeta }P\left( \zeta \right) ,
\label{u-P}
\end{equation}
it is well known that the function $P\left( \zeta \right) $ has the same
limit $p_{0}$ as $\zeta \rightarrow \infty $ along any ray in the region $-%
\frac{3\pi }{2}<\arg \zeta <\frac{3\pi }{2}$.

Any pair of solutions consisting of $u\left( \zeta \right) $ together with
one of the functions $u\left( -\zeta \right) ,u\left( \zeta e^{\pi i}\right)
$ and $u\left( \zeta e^{-\pi i}\right) $ forms a basis in the space of
solutions of (\ref{T-Airy}). So, for example, we have
\begin{equation}
u\left( \zeta e^{\pi i}\right) =Au\left( \zeta e^{-\pi i}\right) +Bu\left(
\zeta \right) ,  \label{MR}
\end{equation}
where $A$ and $B$ are constants. Using (\ref{u-P}) we can rewrite the
relation (\ref{MR}) in the form
\begin{equation}
P\left( \zeta e^{\pi i}\right) =AP\left( \zeta e^{-\pi i}\right)
+Be^{-a\zeta }P\left( \zeta \right) .  \label{MR-Var}
\end{equation}
Noting that all the functions $P\left( \zeta e^{-\pi i}\right) ,P\left(
\zeta \right) $ and $P\left( \zeta e^{\pi i}\right) $ tend to $p_0$ as $%
\zeta \rightarrow +\infty $ and letting $\zeta \rightarrow +\infty $ in (\ref%
{MR-Var}), it follows that $A=1$. Setting $B=T$, we rewrite the relation (%
\ref{MR-Var}) in the following final form
\begin{equation}
P\left( \zeta e^{\pi i}\right) =P\left( \zeta e^{-\pi i}\right) +Te^{-a\zeta
}P\left( \zeta \right) ,  \label{MR-P}
\end{equation}
where
\begin{equation}
T=T\left( a,a_0,a_1,\ldots \right)  \label{T-function}
\end{equation}
is a complex constant which is usually referred to as the \textit{connection
coefficient}, or \textit{Stokes multiplier}. For the particular case of (\ref%
{T-Airy}) given by (\ref{Bessel}) the constant $T$ is known,
\begin{equation}
T=2i\cos \nu \pi .  \label{T-nu}
\end{equation}

In the class of all differential equations with an irregular singular point
at infinity of Poincar\'{e} rank $1$ the pBde should be regarded as an
exceptional or degenerate case in the sense that a pBde generates a pair of
monodromic relations which are not intertwined (each relation involves only
a single function) and the second relation can be derived from the first.
For example, the relation for $P\left( -\zeta \right) $ follows immediately
from (\ref{MR-P}).

Representing $P\left( \zeta \right) $ as a Laplace transform \ \
\begin{equation}
P\left( \zeta \right) =\mathcal{L}\{F\}:=\zeta \int_0^{+\infty }e^{-\zeta
t}F\left( t\right) dt,  \label{L-transform}
\end{equation}
it can be shown using the standard Laplace-Borel technique that $F\left(
t\right) $ satisfies the monodromic relation
\begin{equation}
F\left( \left( t+a\right) e^{\pi i}-a\right) =F\left( \left( t+a\right)
e^{-\pi i}-a\right) -TF\left( -\left( t+a\right) e^{\pi i}\right) .
\label{HMR}
\end{equation}
The relation (\ref{HMR}) can be viewed as dual to (\ref{MR-P}). Moreover, it
can be proved that $F\left( t\right) $ has all the other properties listed
in the second paragraph of Section 1 and also satisfies relation (\ref{***}%
). In particular, for Bessel's case given by (\ref{Bessel}) we have
\begin{equation}
F\left( t\right) ={}\,_2F_1\left( \frac 12-\nu ,\frac 12+\nu ,1;-\frac
t2\right) ,  \label{Bessel's case}
\end{equation}
and as we noted in Section 2 the above properties are well known for this
function. Analysis of the proof of the validity of these properties shows
that instead of the pBde, we actually use the functional equation (\ref{MR-P}%
) and the above properties of $P\left( \zeta \right) $\ which have been
extracted from (\ref{T-Airy}).

Considering the relation (\ref{MR-P}) separately from the differential
equation (\ref{T-Airy}), it turns out that (\ref{MR-P}), while being very
much simpler than (\ref{T-Airy}), yet remains a very rich source of
information.

\subsection{\textbf{The linear space }$S_{a,T}$}

Given $a>0$ and complex $T$ we introduce a linear space $S_{a,T}$ of
functions $P\left( \zeta \right) $ satisfying the conditions:

\begin{enumerate}
\item $P\left( \zeta \right) $ is analytic in the $\zeta $-plane punctured
at the points $\zeta =0$ and $\zeta =\infty ;$

\item Given $P\left( \zeta \right) \in S_{a,T},$ there exists a positive
decreasing function $M_P\left( r\right) $ such that if $\left| \arg \zeta
\right| \leq \pi $\ and $\left| \zeta \right| \geq r$ for any $r,0<r<\infty
, $ then
\begin{equation}
\left| P\left( \zeta \right) \right| \leq M_P\left( r\right) .  \label{one}
\end{equation}

\item Every $P\left( \zeta \right) \in S_{a,T},$ is a solution of the
equation (\ref{MR-P}).
\end{enumerate}

\textbf{Definition 1}. \textit{Given }$P\left( \zeta \right) \in S_{a,T}$%
\textit{\ we call the relation (\ref{MR-P}) the\ Stokes monodromic relation}.

If $P\left( \zeta \right) $ is given by (\ref{T-Airy}) and (\ref{u-P}) and $%
T $ is given by (\ref{T-function}) then $P\left( \zeta \right) \in S_{a,T}.$
This shows that for any pair $\left( a,T\right) $ there exists a non-trivial
$P\left( \zeta \right) \in S_{a,T}$. For instance, given $T$ with $\nu $
defined by (\ref{T-nu}), we represent a solution $u\left( \zeta \right)
=u_{\nu }\left( \zeta \right) $ of the Bessel \ \ equation (\ref{Bessel})
decaying on the positive ray in the form $u_{\nu }\left( \zeta \right)
=e^{-\zeta }P_{\nu }\left( \zeta \right) .$ Clearly $P_{\nu }\left( \frac{a}{%
2}\zeta \right) \in S_{a,T}$.

Despite the fact that Definition 1 was suggested by a differential equation,
in order to distinguish those techniques which depend on the equation itself
from those which do not we have introduced a more general object. The space $%
S_{a,T}$ with independent constants $T$ and $a$ may contain elements not
associated with such a differential equation. Indeed, it is easy to check
that given an entire function $\Phi \left( z\right) $\ and $P\left( \zeta
\right) \in S_{a,T}$\ we have $\Phi \left( 1/\zeta ^{2}\right) P\left( \zeta
\right) \in S_{a,T}$.

It follows from the definition of $S_{a,T}$ that any element $P\left( \zeta
\right) \in S_{a,T}$ remains bounded at infinity in the region $\left| \arg
\zeta \right| \leq \frac{3\pi }2,$ and that in any sectorial region of the $%
\zeta $-plane the function $P\left( \zeta \right) $ is of exponential growth
at infinity with exponent $a$. More precisely, we claim that the following
statements are valid.

Starting with $S_1=1,$ $T_1=T$ we introduce two sequences $\left\{
S_k\right\} _{k=1}^\infty $ and $\left\{ T_k\right\} _{k=1}^\infty $ which
are determined by the following recurrence relations
\begin{align}
& S_{2m+1}=S_{2m},\,\,S_{2m+2}=S_{2m+1}+T_{2m+1}T,  \label{Recurrence-S} \\
& T_{2m+1}=T_{2m}+S_{2m}T,\,\,T_{2m+2}=T_{2m+1}.  \label{Recurrence-T}
\end{align}

\textbf{Lemma 1}. \textit{Let }$P\left( \zeta \right) \in S_{a,T}$\textit{\
and let sequences }$\left\{ S_{k}\right\} _{k=1}^{\infty }$\textit{\ and }$%
\left\{ T_{k}\right\} _{k=1}^{\infty }$\textit{\textit{\ be defined by (\ref%
{Recurrence-S}) and (\ref{Recurrence-T}).Then for all integers }}$m\geq -1$%
\textit{\ the following monodromic relations are valid}
\begin{equation}
P\left( \zeta e^{\left( 2m+1\right) \pi i}\right) =S_{2m+1}P\left( \zeta
e^{-\pi i}\right) +T_{2m+1}e^{-a\zeta }P\left( \zeta \right)  \label{odd}
\end{equation}%
\textit{and}
\begin{equation}
P\left( \zeta e^{\left( 2m+2\right) \pi i}\right) =S_{2m+2}P\left( \zeta
\right) +T_{2m+2}e^{a\zeta }\left( P\left( \zeta e^{-\pi i}\right) \right) .
\label{even}
\end{equation}%
\textit{These relations can be extended to all negative values of }$k$.

\textbf{Illustration. }One can check that for the Bessel case (\ref{Bessel})
the relations (\ref{odd}) and (\ref{even}) with coefficients satisfying (\ref%
{Recurrence-S}) and (\ref{Recurrence-T}) are identical to the formulae for
analytical continuations of Hankel's functions (4.13) and (4.14) from \cite%
{Olver}, respectively.

\textbf{Proof.} We will use mathematical induction on $m$ to prove that the
relations (\ref{odd}) and (\ref{even}) are satisfied for some coefficients $%
S_{2m+1},$ $T_{2m+1},$ $S_{2m+2}$ and $T_{2m+2},$ and in doing so we obtain
the recurrence relations (\ref{odd}) and (\ref{even}).

For $m=-1$ we have the trivial instance of (\ref{odd}), with $S_{-1}=1$ and $%
T_{-1}=0$; and the trivial instance of (\ref{even}), with $S_{0}=1$ and $%
T_{0}=0.$ For $m=0$ (\ref{odd}) follows immediately from (\ref{MR-P}), along
with the initial conditions $S_{1}=1,$ $T_{1}=T.$ From (\ref{MR-P}) we also
have, using (\ref{MR-P}),
\begin{equation*}
P\left( \zeta e^{2\pi i}\right) =P\left( \left( \zeta e^{\pi i}\right)
e^{\pi i}\right) =P\left( \left( \zeta e^{\pi i}\right) e^{-\pi i}\right)
+Te^{-a\zeta e^{\pi i}}P\left( \zeta e^{\pi i}\right) =P\left( \zeta \right)
+Te^{a\zeta }P\left( \zeta e^{\pi i}\right) ,
\end{equation*}
so that
\begin{equation*}
P\left( \zeta e^{2\pi i}\right) =P\left( \zeta \right) +Te^{a\zeta }P\left(
\zeta e^{\pi i}\right) .
\end{equation*}
Now applying (\ref{MR-P}) again gives
\begin{equation*}
P\left( \zeta e^{2\pi i}\right) =P\left( \zeta \right) +Te^{a\zeta }\left(
P\left( \zeta e^{-\pi i}\right) +Te^{-a\zeta }P\left( \zeta \right) \right)
\end{equation*}
which yields finally
\begin{equation*}
P\left( \zeta e^{2\pi i}\right) =\left( 1+T^{2}\right) P\left( \zeta \right)
+Te^{a\zeta }P\left( \zeta e^{-\pi i}\right) .
\end{equation*}
This shows, incidentally, that $S_{2}=1+T^{2}$ and $T_{2}=T,$ and so the
recurrence relations hold when $m=0.$

Suppose now that odd and even hold for some value of $m;$ we can use a
similar process to prove that they also hold for $m+1.$ Thus, using (\ref%
{even}) we have
\begin{equation*}
P\left( \zeta e^{\left( 2\left( m+1\right) +1\right) \pi i}\right)
=S_{2m+2}P\left( \zeta e^{\pi i}\right) +T_{2m+2}e^{a\zeta e^{\pi i}}P\left(
\zeta \right)
\end{equation*}
Applying (\ref{MR-P}) gives
\begin{equation*}
P\left( \zeta e^{\left( 2m+3\right) \pi i}\right) =S_{2m+2}\left( P\left(
\zeta e^{-\pi i}\right) +Te^{-a\zeta }P\left( \zeta \right) \right)
+T_{2m+2}e^{-a\zeta }P\left( \zeta \right)
\end{equation*}
from which we obtain
\begin{equation*}
P\left( \zeta e^{\left( 2m+3\right) \pi i}\right) =S_{2m+2}P\left( \zeta
e^{-\pi i}\right) +\left( T_{2m+2}+S_{2m+2}T\right) e^{-a\zeta }P\left(
\zeta \right) ,
\end{equation*}
and also the corresponding recurrence relations.

In a similar way it can be shown that
\begin{equation*}
P\left( \zeta e^{\left( 2m+4\right) \pi i}\right) =\left(
S_{2m+3}+T_{2m+3}T\right) P\left( \zeta \right) +T_{2m+3}e^{a\zeta }P\left(
\zeta e^{-\pi i}\right) ,
\end{equation*}
and thus that the last pair of recurrence relations hold. $\blacktriangle $

The next statement is an immediate corollary of Lemma 1.

\textbf{Lemma 2. }\textit{Let }$P\left( \zeta \right) \in S_{a,T}$\textit{.
Then for every sector }
\begin{equation*}
S\left( r,\theta \right) :=\left\{ \zeta :|\zeta |\geq r,|\arg \zeta |\leq
\theta \right\}
\end{equation*}
\textit{\ there exists a positive constant }$M_P\left( \theta ,r,T\right) $%
\textit{\ such that for }$\zeta \in S\left( r,\theta \right) $%
\begin{equation}
\left| P\left( \zeta \right) \right| \leq M_P\left( \theta ,r,T\right) \exp
\left( a\left| \zeta \right| \right) .  \label{global}
\end{equation}
\textbf{Proof. }The validity of (\ref{global}) follows from relations (\ref%
{odd}), (\ref{even}) and the inequality (\ref{one}).$\blacktriangle $

We note that every $P\left( \zeta \right) \in S_{a,T}$ can be provided with
an expansion of the form $\sum_{k=0}^\infty p_k/\zeta ^k,$ with complex
coefficients $p_k,$ valid in the region $-\frac{3\pi }2<\arg \zeta <\frac{%
3\pi }2$. Introducing the \textit{remainders}
\begin{equation}
\mathrm{P}_n\left( \zeta \right) :=P\left( \zeta \right) -\sum_{k=0}^{n-1}%
\frac{p_k}{\zeta ^k},n=0,1,\ldots ,  \label{remainders}
\end{equation}
it is possible to prove that

\textbf{Theorem 1}. \textit{For }$P\left( \zeta \right) \in S_{a,T}$ \textit{%
the following estimates are valid for every }$n=0,1,\ldots $\textit{, and
for }$\left\vert \zeta \right\vert \geq r$

\begin{itemize}
\item
\begin{equation}
\left\vert \arg \zeta \right\vert \leq \frac{\pi }{2}\Rightarrow \left\vert
\mathrm{P}_{n}\left( \zeta \right) \right\vert \leq M_{P}\left( r\right)
\frac{n!}{a^{n}\left\vert \zeta \right\vert ^{n}}  \label{errors-(0)}
\end{equation}

\item
\begin{equation}
\left\vert \arg \zeta \right\vert \leq \pi \Rightarrow \left\vert \mathrm{P}%
_{n}\left( \zeta \right) \right\vert \leq M_{P}\left( r\right) \frac{n!\sqrt{%
n+3}}{a^{n}\left\vert \zeta \right\vert ^{n}}{,}  \label{errors-1}
\end{equation}

\item
\begin{equation}
\pi \leq |\arg \zeta |<\frac{3\pi }{2}\Rightarrow \left\vert \mathrm{P}%
_{n}\left( \zeta \right) \right\vert \leq 2M_{P}\left( r\right) \frac{n!%
\sqrt{n+3}}{a^{n}|\Re \{\zeta \}|^{n}},  \label{errors-2}
\end{equation}%
\textit{where }$M_{P}\left( r\right) $\textit{\ is a constant which is
proportional to that from} (\ref{one}).
\end{itemize}

It follows immediately from (\ref{errors-1}) that for $n=0,1,\ldots $%
\begin{equation}
\left\vert p_{n}\right\vert \leq M_{P}\frac{n!}{a^{n}},M_{P}=\inf_{0<r<%
\infty }M_{P}\left( r\right) .  \label{est-p-n}
\end{equation}

The proof of Theorem 1 is based on the techniques of the paper \cite{Gur},
Section 5, using the integral representations (110) and (111) given there
and Theorem 2 of the current paper. We will not prove this in detail here
since a more general result will be presented in our paper \cite{GGM-2}.

The reason for the deterioration when passing from (\ref{errors-1}) to (\ref%
{errors-2}) is the appearance of exponentially small terms upon crossing the
\textit{Stokes rays} $\arg \zeta =\pm \pi .$ This, which is a manifestation
of the Stokes Phenomenon, will be discussed briefly in the next section.

We note finally that Lemma 1 and Lemma 2 allow us to extend the error bounds
for $P\left( \zeta \right) $ given by (\ref{errors-1}) to the whole Riemann
surface $-\infty <\arg \zeta <\infty $ by adding to the remainders
appropriate exponentially small terms upon crossing the Stokes rays, $\arg
\zeta =\pm m\pi, \,m\in \mathbb{N}$.

\subsection{\textbf{Exponentially small terms and Stokes' Phenomenon}}

The relation (\ref{MR-P}) provides an immediate clarification of the r\^{o}%
le of exponentially small terms in Stokes phenomenon, as we will now show.

Comparison of relations (\ref{errors-1}) and (\ref{errors-2}) shows that
upon crossing the rays $\arg \zeta =\pm \pi $ the approximation given by (%
\ref{errors-1}) begins to deteriorate and the greater the deviation from the
ray the greater this deterioration becomes.

The cause of the above deterioration is hidden in relation (\ref{MR-P}).
Indeed, relation (\ref{MR-P}) generates the system of relations

\begin{equation}
\mathrm{P}_n\left( \zeta e^{\pi i}\right) =\mathrm{P}_n\left( \zeta e^{-\pi
i}\right) +Te^{-a\zeta }P\left( \zeta \right) {,}  \label{n-remainders}
\end{equation}
where $n=0,1,\ldots $ and $\mathrm{P}_0\left( \zeta \right) =P\left( \zeta
\right) $. Let us introduce functions
\begin{align}
E^{+}\left( \zeta \right) & =Te^{a\zeta }P\left( \zeta e^{-\pi i}\right)
,\pi \leq \arg \zeta <\frac{3\pi }2  \label{E+small} \\
E^{-}\left( \zeta \right) & =-Te^{a\zeta }P\left( \zeta e^{\pi i}\right) ,-%
\frac{3\pi }2<\arg \zeta \leq -\pi {.}  \label{E-small}
\end{align}
These functions decay exponentially in the regions $\left\{ \zeta :\pi \leq
\arg \zeta <\frac{3\pi }2\right\} $\ and $\left\{ \zeta :-\frac{3\pi }2<\arg
\zeta \leq -\pi \right\} ,$ respectively, as $\zeta \rightarrow \infty $.
Using (\ref{n-remainders}) we observe that subtracting these exponentially
small terms from the remainders, we retrieve the accuracy of (\ref{errors-1}%
):
\begin{align}
\pi & \leq \arg \zeta \leq \frac{3\pi }2\Longrightarrow \left| \mathrm{P}%
_n\left( \zeta \right) -E^{+}\left( \zeta \right) \right| \leq M_P\left(
r\right) \frac{\sqrt{n+3}n!}{a^n|\zeta ^n|},  \label{errors-(+)} \\
-\frac{3\pi }2& \leq \arg \zeta \leq -\pi \Longrightarrow \left| \mathrm{P}%
_n\left( \zeta \right) -E^{-}\left( \zeta \right) \right| \leq M_P\left(
r\right) \frac{\sqrt{n+3}n!}{a^n|\zeta ^n|},  \label{errors-(-)}
\end{align}

Indeed, it follows from (\ref{n-remainders}), (\ref{E+small}) and (\ref%
{E-small}) that for $n=0,1,\ldots $%
\begin{align*}
{\pi }& {\leq \arg \zeta \leq \frac{3\pi }2}\Longrightarrow \arg \left(
\zeta e^{-2\pi i}\right) \in \left[ -\pi ,\pi \right] :\mathrm{P}_n\left(
\zeta \right) -E^{+}\left( \zeta \right) =\mathrm{P}_n\left( \zeta e^{-2\pi
i}\right) , \\
{-\frac{3\pi }2}& {<\arg \zeta \leq -\pi \Longrightarrow \arg \left( \zeta
e^{2\pi i}\right) \in \left[ -\pi ,\pi \right] :\mathrm{P}_n\left( \zeta
\right) -E^{-}\left( \zeta \right) =\mathrm{P}_n\left( \zeta e^{2\pi
i}\right) ,}
\end{align*}
which, using (\ref{errors-1}) yields (\ref{errors-(+)}) and (\ref{errors-(-)}%
).

Note that if the monodromic relation is generated by the differential
equation (\ref{T-Airy}) then evaluation of the exponentially small terms
reduces to evaluation of the constant $T.$ To clarify the cause for the
above phenomenon we must turn to the dual complex plane and to consider a
dual space $H_{a,T}$ of functions $F\left( t\right) $ analytic in the dual
plane.

\subsection{\textbf{The linear space }$H_{a,T}$}

Given $a>0$ and complex $T$ we introduce a second linear space $H_{a,T}$ of
functions $F\left( t\right) $ satisfying the conditions:

\begin{itemize}
\item (i) $F\left( t\right) $\ is an analytic function in the complex $t$%
-plane punctured at the two finite points $t=0$\ and $t=-a,$ and with
exponential growth of minimal type at $\infty $ in any sectorial region $%
\{t:-\infty <\alpha <\arg t<\beta <+\infty \};$

\item (ii) there exists a branch of $F(t)$\ which is analytic in the $t$%
-plane cut along $\left( -\infty ,-a\right] ;$

\item (iii) the branch given by (ii) satisfies the estimate $F(t)=O\left(
\log \left( t+a\right) \right) $ as $t\rightarrow -a,\left| t+a\right|
<a,\left| \arg \left( t+a\right) \right| <\pi $;

\item (iv) the branch given by (ii) satisfies the monodromic relation (\ref%
{HMR}) given by
\begin{equation}
F\left( \left( t+a\right) e^{\pi i}-a\right) =F\left( \left( t+a\right)
e^{-\pi i}-a\right) -TF\left( -\left( t+a\right) e^{\pi i}\right) ,
\label{HMR-*}
\end{equation}
where $-a<t<+\infty $ , and $T$ is the constant from (\ref{MR-P}).

Note that the points $\left( t+a\right) e^{\pi i}-a$ and $\left( t+a\right)
e^{-\pi i}-a$ belong to the upper and lower banks of the cut, respectively,
and the relation (\ref{HMR-*}) can be rewritten as
\begin{equation*}
F\left( te^{\pi i}-a\right) =F\left( te^{-\pi i}-a\right) -TF\left( -te^{\pi
i}\right) ,t>0,
\end{equation*}
as well as in the equivalent form
\begin{equation}
F\left( \left( t+a\right) e^{-2\pi i}-a\right) =F\left( t\right) +TF\left(
-\left( t+a\right) \right) ,  \label{HMR+}
\end{equation}
where $t$ belongs to the upper bank of the cut.
\end{itemize}

\textbf{Definition 2}. \textit{Given }$F\left( t\right) \in H_{a,T}$\textit{%
\ we call} \textit{the relation (\ref{HMR-*})\ the hypergeometric monodromic
relation.}

\textbf{Remark 1}. \textit{Since function }$F\left( t\right) $\textit{\ is
multi-valued, it is necessary to interpret the expression $F\left(
-t\right)\, $ ``carefully''. If $F_0\left( t\right) $ is a branch of $%
F\left( t\right) $ which is given by (ii) of Definition 2, then $F_0\left(
-t\right) $ is a single-valued analytic function in the $t$-plane cut along $%
\left[ a,+\infty \right) $ such that}
\begin{equation*}
F_0(-t)=%
\begin{cases}
F_0\left(t{e^{\pi i}}\right ) & \text{if } \Im {t}<0\,, \\
F_0\left( t{e^{-\pi i}}\right) & \text{if } \Im {t}>0\,.%
\end{cases}%
\end{equation*}
\textit{We define }$F\left( -t\right) $\textit{\ as an analytical
continuation of }$F_0\left( -t\right) $\textit{\ to the complex }$t$\textit{%
-plane punctured at the two points }$t=0$\textit{\ and }$t=-a$\textit{. The
function }$F\left( -t-a\right) $\textit{\ satisfies condition (i).}

Below we will also use $F\left( t\right) $ to denote the branch given by
condition (ii), where this will cause no confusion.

For the case given by (\ref{Bessel's case}) setting $F_{\nu }\left( t\right)
=\, _{2}F_{1}\left( \frac{1}{2}-\nu ,\frac{1}{2}+\nu ,1;-\frac{t}{a}\right) $%
, where $a=\frac{1}{2}-\nu$ and $T=2i\cos \nu \pi, $ one can check that $%
F_{\nu }\left( t\right) \in H_{a,T}.$

In the next section we will show that there is a one-to-one correspondence
between the linear spaces $H_{a,T}$ and $S_{a,T}$.

\section{\textbf{The Duality theorem}}

Defining the Laplace transform operator $\mathcal{L}$ by (\ref{L-transform})
we claim that

\textbf{Theorem 2}. \textit{Given }$a$\textit{\ and }$T,$ \textit{the
operator }$\mathcal{L}:H_{a,T}\rightarrow S_{a,T}$\textit{\ is a bijection}.

We present the proof in two parts.

\subsection{\textbf{Part I: }$\mathcal{L}H_{a,T}\subset S_{a,T}$}

\textbf{Proof. }Let $F\left( t\right) \in H_{a,T}$ and let $P\left( \zeta
\right) $ be given by (\ref{L-transform}). Then it follows from conditions
(i) and (ii) of section 2.3 that $P\left( \zeta \right) $ admits an
analytical continuation from the right half-plane to the $\zeta $-plane cut
along the interval $\left( -\infty ,0\right) $ as a function bounded at $%
\infty $ and continuous in the cut plane except for $\zeta =0$. Moreover,
for $-\pi \leq \theta \leq \pi $ and for $0<\rho <+\infty ,$ using Cauchy's
theorem, the integral representation
\begin{equation}
P\left( \rho e^{i\theta }\right) =\rho e^{i\theta }\int_0^{\infty \cdot
e^{-i\theta }}e^{-\rho e^{i\theta }t}F\left( t\right) dt  \label{P-F-theta}
\end{equation}
is valid, where the integral is absolutely convergent. Setting $\theta =\pm
\pi $ in (\ref{P-F-theta}), we retain the previous convergence due to the
condition (iii) of 2.3. Subtracting the second integral from the first we
have
\begin{equation}
P\left( \rho e^{\pi i}\right) -P\left( \rho e^{-\pi i}\right) =-\rho
\int_0^{\infty \cdot e^{-\pi i}}e^{\rho t}F\left( t\right) dt+\rho
\int_0^{\infty \cdot e^{\pi i}}e^{\rho t}F\left( t\right) dt.
\label{the first}
\end{equation}

Due to condition (ii) $F\left( t\right) $ is analytic in the circle $%
\left\vert t\right\vert <a,$ and\textbf{\ }thus for $t\in \left( -a,0\right)
$ we have $F\left( t\right) =F\left( t^{\ast }\right) $, where $t^{\ast
}=\left( t+a\right) e^{-2\pi i}-a.$ Therefore, writing $\int_{-a}^{-\infty }$%
\ for an integral along the upper bank of the cut, (\ref{the first}) can be
rewritten in the form
\begin{align}
& P\left( \rho e^{\pi i}\right) -P\left( \rho e^{-\pi i}\right) =-\rho
\int_{-a}^{-\infty }e^{\rho t}F\left( t^{\ast }\right) dt+\rho
\int_{-a}^{-\infty }e^{\rho t}F\left( t\right) dt  \label{the second 1} \\
=-& \rho \int_{-a}^{-\infty }e^{\rho t}\left( F\left( t^{\ast }\right)
-F\left( t\right) \right) dt.  \label{the second}
\end{align}
Using the monodromic relation (\ref{HMR-*}) in the form (\ref{HMR+}), the
last integral in (\ref{the second}) can be rewritten as
\begin{equation}
-T\rho \int_{-a}^{-\infty }e^{\rho t}F\left( -\left( t+a\right) \right)
dt=T\rho e^{-a\rho }\int_{-\infty }^{-a}e^{\rho \left( t+a\right) }F\left(
-\left( t+a\right) \right) dt.  \label{f-1}
\end{equation}
Finally we have
\begin{equation}
T\rho e^{-a\rho }\int_{-\infty }^{0}e^{\rho t}F\left( -t\right) dt=T\rho
e^{-a\rho }\int_{0}^{+\infty }e^{-\rho t}F\left( t\right) dt=Te^{-a\rho
}P\left( \rho \right) .  \label{f-2}
\end{equation}
The chain of relations (\ref{the second})--(\ref{f-2}) shows that $P\left(
\zeta \right) $ satisfies the condition (iii) of definition 1, and thus, $%
P\left( \zeta \right) \in S_{a,T}.\blacktriangle $

\subsection{\textbf{II.} $S_{a,T}\subset \mathcal{L}H_{a,T}$}

We assume now that $P\left( \zeta \right) \in S_{a,T}.$ We set $F\left(
t\right) $ equal to the Borel transform of $P\left( \zeta \right) ,$ so that
$P\left( \zeta \right) =\mathcal{L}\left\{ F\right\} ,$ and prove, in five
steps, that $F\left( t\right) \in H_{a,T}.$

\textbf{(1)} \textbf{Integral representation for $F\left( t\right) $%
\thinspace for $\left| \arg t\right| \mathbf{<}\frac \pi 2$.} Given $r>0$
and $\theta >0,$ we introduce the contour $\gamma _\theta \left( r\right) $
with an anti-clockwise orientation as follows:
\begin{equation}
\gamma _\theta \left( r\right) =l_{-\theta }\left( r\right) \cup C_\theta
\left( r\right) \cup l_\theta \left( r\right) ,  \label{gamma-ro}
\end{equation}
where
\begin{equation}
l_\theta \left( r\right) =\left\{ \zeta :\zeta =\rho e^{i\theta },r<\rho
<\infty \right\}  \label{theta-ro}
\end{equation}
and
\begin{equation}
C_\theta \left( r\right) =\left\{ \zeta :-\theta <\arg \zeta <\theta ,\left|
\zeta \right| =r\right\} .  \label{C-r}
\end{equation}
For $t>0$ the function $F\left( t\right) $ can then be represented in the
form
\begin{equation}
F\left( t\right) =\frac 1{2\pi i}\int_{\gamma _\theta \left( r\right)
}e^{t\zeta }P\left( \zeta \right) \frac{d\zeta }\zeta ,  \label{Borel}
\end{equation}
with any $\frac \pi 2<\theta <\frac{3\pi }2$ and $r>0.$

Below we set $\theta =\pi $ and omit the subscript $\pi $ in $C_\pi \left(
r\right) $ and in $\gamma _\pi \left( r\right) $ if this does not lead to
confusion.

First we prove that this integral converges for $t>0.$ Then we demonstrate
how, by changing the integral representation, we can obtain an analytical
continuation of $F\left( t\right) $ that allows us to check the validity of
conditions (i)-(iii) of 2.3. Finally we show that the function $F\left(
t\right) $ satisfies the relation (\ref{HMR-*}). Setting
\begin{equation}
F^{+}\left( t,r\right) =\frac 1{2\pi i}\int_{l_\pi \left( r\right)
}e^{t\zeta }P\left( \zeta \right) \frac{d\zeta }\zeta ,  \label{F+}
\end{equation}
\begin{equation}
F^{-}\left( t,r\right) =\frac 1{2\pi i}\int_{l_{-\pi }\left( r\right)
}e^{t\zeta }P\left( \zeta \right) \frac{d\zeta }\zeta  \label{F-}
\end{equation}
and
\begin{equation}
F_0\left( t,r\right) =\frac 1{2\pi i}\int_{C\left( r\right) }e^{t\zeta
}P\left( \zeta \right) \frac{d\zeta }\zeta  \label{F0}
\end{equation}
we have
\begin{equation}
F\left( t\right) =F^{+}\left( t,r\right) -F^{-}\left( t,r\right) +F_0\left(
t,r\right) .  \label{Fsum}
\end{equation}
Observe, that $F_0\left( t,r\right) $ admits an analytical continuation from
the positive ray to the whole $t$-plane as an entire function of $t$. Since $%
P\left( \zeta \right) \in S_{a,T},$ the relation (\ref{one}) yields the
following estimate valid in the whole $t$-plane
\begin{equation}
\left| F_0\left( t,r\right) \right| \leq M_P\left( r\right) e^{r\left|
t\right| }  \label{F-0-est}
\end{equation}
meaning that $F_0\left( t,r\right) $ is an entire function of exponential
type $r.$ On the other hand, for the functions $F^{\pm }\left( t,r\right) $
given by (\ref{F+}) and (\ref{F-}) the inequality (\ref{one}) yields the
estimates
\begin{equation}
\left| F^{\pm }\left( t,r\right) \right| \leq \frac{M_P\left( r\right) }{%
2\pi r}\int_r^{+\infty }e^{-\zeta t}d\zeta =\frac{M_P\left( r\right) }{2\pi r%
}e^{-rt},0<t<+\infty .  \label{F+-Estimate}
\end{equation}
The estimates (\ref{F-0-est}) and (\ref{F+-Estimate}) show that integrals (%
\ref{F+}) -- (\ref{F0}) are absolutely convergent, so that the function $%
F\left( t\right) $ is well defined for positive $t,$ and that its growth on
the positive ray is determined essentially by estimate (\ref{F-0-est}). It
follows from (\ref{F-0-est}) and the preceding comments that the integral in
(\ref{L-transform}) is absolutely convergent for $\Re \left\{ \zeta \right\}
>r.$ Since $r$ is an arbitrarily small positive number it follows that this
integral is absolutely convergent for all $\zeta $ such that $\left| \arg
\zeta \right| <\frac \pi 2.$ Therefore formula (\ref{Borel}) can be
considered as the converse for the Laplace transform (\ref{L-transform}),
which proves the validity of the representation
\begin{equation}
P\left( \zeta \right) =\zeta \int_0^\infty e^{-s\zeta }F\left( s\right)
ds,-\frac \pi 2<\arg \zeta <\frac \pi 2.  \label{Laplace Transform}
\end{equation}

\textbf{(2) Representation for }$F\left( t\right) $\textbf{\ in the region }$%
\left\{ t:\left| t\right| >a\right\} .$\textit{\ }Let $\theta $ be any
angle, $-\infty <\theta <\infty $. To obtain the desired analytical
continuation of $F\left( t\right) $ we rotate the path of integration in (%
\ref{Borel}), with $\theta =\pi $, about the origin through an angle $%
-\theta $ while simultaneously increasing $\arg t$ by $\theta $. Lemma 1
yields
\begin{equation}
\left| P\left( \zeta \right) \right| \leq M\left( \theta ,r,T\right) \exp
\left( a\left| \zeta \right| \right) ,  \label{P-growth}
\end{equation}
where the positive parameter $M\left( \theta ,r,T\right) $ depends on $%
\theta $, $r$ and $T$. Following the above rotation the integral in the form
\begin{equation}
F\left( te^{i\theta }\right) :=\frac 1{2\pi i}\int_{\gamma \left( r\right)
\cdot e^{-i\theta }}e^{te^{i\theta }\zeta }P\left( \zeta \right) \frac{%
d\zeta }\zeta ,0<t<\infty ,  \label{F - large t rep}
\end{equation}
retains the absolute convergence of the integral (\ref{Borel}). Thus, using
Cauchy's theorem, it can be shown that (\ref{F - large t rep}) provides an
analytical continuation of (\ref{Borel}). Since $\theta $ is arbitrary, this
shows that $F\left( t\right) $ admits an analytical continuation to the
exterior of the circle $\left\{ t:\left| t\right| =a\right\} $, and (\ref{F
- large t rep}) provides an integral representation for $F\left( t\right) $
in the region $\left\{ t:\left| \arg t\right| <\pi \right\} $ and also in
the region $\left\{ t:\left| t\right| >a,-\infty <\arg t<\infty \right\} .$
It also follows from (\ref{F - large t rep}) and (\ref{P-growth}) that given
$0<r<\infty $ in any sectorial region $\left\{ t:|t|\geq R>a,\alpha <\arg
t<\beta \right\} ,$\ of the $t$-plane the following estimate is valid
\begin{equation}
\left| F\left( t\right) \right| \leq Ke^{r\left| t\right| },
\label{minimal growth}
\end{equation}
where the positive parameter $K$\ depends on $R,\alpha ,\beta ,$ and $r.$

\textbf{(3)} \textbf{Representation of }$F\left( t\right) $\textbf{\ in the
region }$\left( t:\left\vert \arg \left( t+a\right) \right\vert <\frac{\pi }{%
2}\right) $. Next we prove that $F\left( t\right) $ is analytic in the
half-plane $\Re z>-a$ and therefore it satisfies the conditions (i) and (ii)
of 2.3. In order to prove this we first assume that $0<t<\infty $ and derive
another integral representation of $F\left( t\right) $. Now, returning to
the representation (\ref{Borel}) and making, in (\ref{F+}) and (\ref{F-}),
the change of variable suggested by (\ref{theta-ro}) for $\theta =\pm \pi $,
it follows that (\ref{Borel}) can be written in the form

\begin{equation}
F\left( t\right) =F_{0}\left( t,r\right) +\frac{1}{2\pi i}\int_{r}^{+\infty
}e^{-t\zeta }\left( P\left( \zeta e^{\pi i}\right) -P\left( \zeta e^{-\pi
i}\right) \right) \frac{d\zeta }{\zeta },  \label{on the cut}
\end{equation}
where\ $F_{0}\left( t,r\right) $ is given by (\ref{F0}).

Using the relation
\begin{equation*}
P\left( \zeta e^{i\pi }\right) -P\left( \zeta e^{-i\pi }\right) =Te^{-a\zeta
}P\left( \zeta \right),
\end{equation*}
which follows from (\ref{MR-P}), we can represent $F\left( t\right) ,\,t>0$
in the form
\begin{equation}
F\left( t\right) =F_{0}\left( t,r\right) +\frac{T}{2\pi i}\int_{r}^{+\infty
}e^{-\left( t+a\right) \zeta }P\left( \zeta \right) \frac{d\zeta }{\zeta }.
\label{SBorel}
\end{equation}
Formula (\ref{SBorel}) shows that $F\left( t\right) $ can be continued
analytically from the positive ray to the half-plane $\left\vert \arg \left(
t+a\right) \right\vert <\frac{\pi }{2}$. This completes the proof of step
\textbf{(3).}

Steps \textbf{(1)-(3) }prove that $F\left( t\right) $ is analytic in the
whole $t$-plane punctured at $t=0$ and $t=-a,$ and formula (\ref{SBorel})
represents a branch of $F\left( t\right) $ which is analytic in the plane
cut along the interval $\left( -\infty ,-a\right] $. We keep the notation $%
F\left( t\right) $ for this branch, if it does not lead to confusion.
Clearly the function $F\left( t\right) $ is analytic in the circle $%
\left\vert t\right\vert <a$. Applying Watson's lemma to the relation (\ref%
{Laplace Transform}) we can provide every element $P\left( \zeta \right) \in
S_{a,T}$ with an asymptotic expansion of the form $\sum_{k=0}^{\infty
}p_{k}/\zeta ^{k}$, where
\begin{equation}
p_{k}=F^{\left( k\right) }\left( 0\right) ,k=0,1,\ldots ,F=\mathcal{L}^{-1}P.
\label{coef-k}
\end{equation}%
Analytical properties of $F\left( t\right) $ given by (1)-(3) show that the
relation
\begin{equation}
P\left( \zeta \right) -\sum_{k=0}^{n-1}p_{k}/\zeta ^{k}=O\left( \frac{1}{%
\zeta ^{n}}\right) ,\zeta \rightarrow \infty ,  \label{Poincare-asymp}
\end{equation}%
is valid for every $n\in \mathbb{N}$ and for $\zeta $ satisfying $-\frac{\pi
}{2}<\arg \zeta <\frac{\pi }{2}$, and can be extended to every sub-sector of
the region $-\frac{3\pi }{2}<\arg \zeta <\frac{3\pi }{2}.$ Next we use this
fact to prove that $F\left( t\right) $ satisfies condition (iii) of 2.3.

\textbf{(4) Behavior of }$F\left( t\right) $\textbf{\ at }$t=-a$\textbf{. }%
We prove the following\textbf{\ }statement.

\textbf{Lemma 2. }\textit{Let }$P\left( \zeta \right) \in S_{a,T}$\textit{\
and }$F\left( t\right) $\textit{\ be given by (\ref{SBorel}), then there
exists a constant }$A_{0}$\textit{\ such that in the region }$\left\vert
t+a\right\vert <a,\left\vert \arg \left( t+a\right) \right\vert <\pi $%
\textit{\ the following formula is valid}
\begin{equation}
F\left( t\right) =-\frac{Tp_{0}}{2\pi i}\log \left( t+a\right)
+A_{0}+o\left( 1\right) ,t\rightarrow -a.  \label{lemma}
\end{equation}

\textbf{Proof.} To prove (\ref{lemma}) we use (\ref{SBorel}) to estimate $%
F\left( t\right) .$ Note that $F_{0}\left( t,r\right) $ is an entire
function in the $t$-plane and from (\ref{F0}) we may write
\begin{equation*}
F_{0}\left( t,r\right) =\frac{1}{2\pi i}\int\limits_{C\left( r\right)
}e^{\left( t+a\right) \zeta }e^{-a\zeta }P\left( \zeta \right) \frac{d\zeta
}{\zeta }.
\end{equation*}
Since $e^{\left( t+a\right) \zeta }=1+O\left( t+a\right) $ as $t\rightarrow
-a$ it then follows that
\begin{equation}
F_{0}\left( t,r\right) =\frac{1}{2\pi i}\int\limits_{C\left( r\right)
}e^{-a\zeta }P\left( \zeta \right) \frac{d\zeta }{\zeta }+O\left( t+a\right)
,t\rightarrow -a.  \label{F-0-r}
\end{equation}

Let us denote the integral on the right hand side of (\ref{SBorel}) by $%
I\left( t,r\right) $ and observe that it can be written in the form
\begin{equation}
I\left( t,r\right) =\frac T{2\pi i}\int_r^{+\infty }e^{-\left( t+a\right)
\zeta }\frac{P\left( \zeta \right) -p_0}\zeta d\zeta +\frac{Tp_0}{2\pi i}%
\int_r^{+\infty }e^{-\left( t+a\right) \zeta }\frac{d\zeta }\zeta ,
\label{second int}
\end{equation}
where $p_0$ is given by (\ref{coef-k}) with $k=0.$ Using the relation (\ref%
{Poincare-asymp}) with $n=1$ together with $e^{\left( t+a\right) \zeta
}=1+O\left( t+a\right) $ as $t\rightarrow -a$, the first integral in (\ref%
{second int}) can be represented as
\begin{equation}
\frac T{2\pi i}\int_r^{+\infty }e^{-\left( t+a\right) \zeta }\frac{P\left(
\zeta \right) -p_0}\zeta d\zeta =\frac T{2\pi i}\int_r^{+\infty }\frac{%
P\left( \zeta \right) -p_0}\zeta d\zeta +O\left( t+a\right) ,t\rightarrow -a,
\label{estimate}
\end{equation}
where the integral on the right hand-side is absolutely convergent.

Finally we evaluate the second integral in (\ref{second int}) using the
exponential integral
\begin{equation}
E_{1}\left( z\right) =\int_{1}^{\infty }e^{-z\tau }\frac{d\tau }{\tau }
\label{E-1}
\end{equation}
\ and its power series expansion, see formula (\textbf{5.1.11}) in \cite%
{Abramowitz},
\begin{equation}
E_{1}\left( z\right) =-\log z-\gamma -\sum_{m=1}^{\infty }\frac{\left(
-z\right) ^{m}}{\left( m+1\right) m!},  \label{E-1-series}
\end{equation}
where $\left\vert \arg z\right\vert <\pi $ and $\gamma $ is Euler's
constant. Thus, replacing $z$ in (\ref{E-1}) by $r\left( t+a\right) $, and
applying (\ref{E-1-series}) we have, for $\left\vert \arg \left( t+a\right)
\right\vert <\pi, $%
\begin{align}
&\frac{Tp_{0}}{2\pi i}\int_{r}^{+\infty }e^{-\left( t+a\right) \zeta }\frac{%
d\zeta }{\zeta }=\frac{Tp_{0}}{2\pi i}E_{1}\left( r\left( t+a\right) \right)
\notag \\
=&\frac{Tp_{0}}{2\pi i}\left( -\log \left( r\left( t+a\right) \right)
-\gamma -\sum_{m=1}^{\infty }\frac{\left( -r\left( t+a\right) \right) ^{m}}{%
\left( m+1\right) m!}\right) .  \label{using E-1}
\end{align}
Combining (\ref{F-0-r}), (\ref{estimate}) and (\ref{using E-1}), we obtain
the following estimate for $F\left( t\right) $
\begin{equation}
F\left( t\right) =-\frac{TF\left( 0\right) }{2\pi i}\log \left( t+a\right)
+A_{0}+O\left( t+a\right) ,t\rightarrow -a,\left\vert \arg \left( t+a\right)
\right\vert <\pi  \label{at-a}
\end{equation}
where
\begin{equation}
A_{0}=\frac{1}{2\pi i}\int\limits_{C\left( r\right) }e^{-a\zeta }P\left(
\zeta \right) \frac{d\zeta }{\zeta }+\frac{T}{2\pi i}\int_{r}^{+\infty }%
\frac{P\left( \zeta \right) -p_{0}}{\zeta }d\zeta -\frac{Tp_{0}}{2\pi i}%
\left( \gamma +\log r\right) ,  \label{A-0}
\end{equation}
and $r\in \left( 0,\infty \right) .$ While the individual terms on the
righthand side of (\ref{A-0}) are clearly dependent on $r$ it follows from (%
\ref{at-a}), taking the limit as $t\rightarrow -a,$that the coefficient $%
A_{0}$ does not depend on $r$. $\blacktriangle $

It remains to be proven that $F\left( t\right) $ satisfies the monodromic
relation (\ref{HMR-*}).

\textbf{(5) Derivation of the monodromic relation for }$F\left( t\right) .$
The properties of $F\left( t\right) $ given by (i) and (ii) of the section
2.3. which have been proven in the preceding paragraph allow us to use the
relations (\ref{the first}) and (\ref{the second}) and to show that for $%
0<\rho <+\infty $%
\begin{equation*}
P\left( \rho e^{\pi i}\right) -P\left( \rho e^{-\pi i}\right) =\rho
\int_{-a}^{-\infty }e^{\rho t}\left( F\left( t^{\ast }\right) -F\left(
t\right) \right) dt,
\end{equation*}
where $t^{\ast }=\left( t+a\right) e^{-2\pi i}-a.$ On the other hand using
elementary transformations given by (\ref{f-1}) and (\ref{f-2}) we have
\begin{equation*}
Te^{-\rho t}P\left( \rho \right) =T\rho \int_{-a}^{-\infty }e^{\rho
t}F\left( -\left( t+a\right) \right) dt.
\end{equation*}
Since the $P\left( \rho e^{\pi i}\right) -P\left( \rho e^{-\pi i}\right)
=Te^{-\rho t}P\left( \rho \right) $, it follows that $F\left( t^{\ast
}\right) -F\left( t\right) =TF\left( -\left( t+a\right) \right) .$

Steps (\textbf{1})-(\textbf{5}) justify the inclusion $S_{a,T}\subset
\mathcal{L}H_{a,T}$. $\blacktriangle $

In the next section we use the Duality theorem to prove a version of the
linear transformation formula for $F\left( t\right) \in H_{a,T}$.

\section{Linear transformation formula}

The following statements is valid.

\textbf{Theorem 3. }\textit{Assume that }$F\left( t\right) \in H_{a,T}.$
\textit{Then for the branch of }$F\left( t\right) $\textit{\ given by
condition (ii) of 2.3 in the region} $\left| t+a\right| <a,\left| \arg
\left( t+a\right) \right| <\pi $ \textit{the following formula is valid}
\begin{equation}
F\left( t\right) =\sum_{k=0}^\infty A_k\left( t+a\right) ^k-\frac T{2\pi
i}F\left( -\left( t+a\right) \right) \log \left( t+a\right) ,  \label{at -2}
\end{equation}
\textit{where the series in (\ref{at -2}) is absolutely convergent, }$A_0$%
\textit{\ is given by (\ref{A-0}), where }$p_0=F\left( 0\right) $\textit{,
and }$A_k,$\textit{\ }$k\in \mathbb{N},$ \textit{are complex coefficients
which can be found explicitly}.

\textbf{Remark 2}. \textit{For Bessel's case given by (\ref{Bessel})
relation (\ref{at -2}) follows immediately from a degenerate case of the
Euler transformation formula for the hypergeometric function, see }15.3.10%
\textit{\ of \cite{Abramowitz} and Appendix 1. Thus, the relation (\ref{at
-2}) can be viewed as a generalization for the elements of }$H_{a,T}$
\textit{of the linear transformation formula (\ref{log-*}) with $\mathfrak{a}%
+\mathfrak{b}=1.$}

\textbf{Proof of Theorem 3}. The relation (\ref{lemma}) of Lemma 2 can be
rewritten as
\begin{equation}
F\left( t\right) =-\frac{TF\left( 0\right) }{2\pi i}\log \left( t+a\right)
+A_0+o\left( 1\right) ,t\rightarrow -a,  \label{Lemma}
\end{equation}
where $A_0$ is given by (\ref{A-0}). Since $F\left( t+a\right) $ is analytic
in the circle $\left| t+a\right| <a$ we have the relation
\begin{equation*}
F\left( -\left( t+a\right) \right) \log \left( t+a\right) =\left( F\left(
0\right) +O\left( t+a\right) \right) \log \left( t+a\right) \text{ as }%
t\rightarrow -a,
\end{equation*}
which yields
\begin{equation*}
F\left( 0\right) \log \left( t+a\right) =F\left( -\left( t+a\right) \right)
\log \left( t+a\right) +o\left( 1\right) .
\end{equation*}
Substituting the last relation into (\ref{Lemma}) allows us to rewrite it as
\begin{equation}
F\left( t\right) =-\frac T{2\pi i}F\left( -\left( t+a\right) \right) \log
\left( t+a\right) +A_0+o\left( 1\right) ,t\rightarrow -a,  \label{at--a}
\end{equation}
Let us introduce a function $\Phi \left( t\right) $ such that

\begin{equation}
F\left( t\right) =-\frac T{2\pi i}F\left( -\left( t+a\right) \right) \log
\left( t+a\right) +A_0+\Phi \left( t+a\right) ,t\rightarrow -a.
\label{with phy}
\end{equation}
Clearly $\Phi \left( t\right) $ is analytic in the $t$-plane punctured at
the two points $t=0$ and $t=a$. On the other hand the relation (\ref{at--a})
shows that\
\begin{equation}
\Phi \left( t+a\right) =o\left( 1\right) ,t\rightarrow -a.  \label{o(1)}
\end{equation}

Since $\left\vert t+a\right\vert <a,\left\vert \arg \left( t+a\right)
\right\vert <\pi $ we can continue analytically both sides of (\ref{with phy}%
) by rotating $t+a$ through angles $\pi $ and $-\pi $ about the origin. This
allows us to rewrite (\ref{with phy}) as
\begin{align}
& F\left( \left( t+a\right) e^{\pi i}-a\right) ={A_{0}}+\Phi \left( \left(
t+a\right) e^{\pi i}\right)  \notag \\
+& \frac{T}{2\pi i}F\left( -\left( t+a\right) e^{\pi i}\right) \log \left(
\left( t+a\right) e^{\pi i}\right) ,  \label{with phy +}
\end{align}
and also as
\begin{align}
& F\left( \left( t+a\right) e^{-\pi i}-a\right) =A_{0}+\Phi \left( \left(
t+a\right) e^{-\pi i}\right)  \notag \\
-& \frac{T}{2\pi i}F\left( -\left( t+a\right) e^{-\pi i}\right) \log \left(
\left( t+a\right) e^{-\pi i}\right) .  \label{with phy -}
\end{align}
Subtracting (\ref{with phy -})from (\ref{with phy +}) and using the
monodromic relation (\ref{HMR+}) gives
\begin{equation*}
F\left( \left( t+a\right) e^{\pi i}-a\right) -F\left( \left( t+a\right)
e^{-\pi i}-a\right) =-TF\left( -\left( t+a\right) e^{\pi i}\right) .
\end{equation*}

A straightforward calculation using
\begin{equation*}
F\left( -\left( t+a\right) e^{\pi i}\right) \equiv F\left( -\left(
t+a\right) e^{-\pi i}\right) ,\left\vert t+a\right\vert <a,
\end{equation*}
shows that
\begin{equation*}
-TF\left( -\left( t+a\right) e^{\pi i}\right) =-\frac{T}{2\pi i}F\left(
-\left( t+a\right) e^{\pi i}\right) \log \left( e^{2\pi i}\right) +\Phi
\left( \left( t+a\right) e^{\pi i}\right) -\Phi \left( \left( t+a\right)
e^{-\pi i}\right) .
\end{equation*}
Canceling equal terms, the last relation can be simplified and we have
\begin{equation*}
\Phi \left( \left( t+a\right) e^{\pi i}\right) =\Phi \left( \left(
t+a\right) e^{-\pi i}\right)
\end{equation*}
which, using (\ref{o(1)}), shows that $\Phi \left( t\right) $ is analytic
and single-valued in the circle $\left\vert t+a\right\vert <a$ and that $%
\Phi \left( 0\right) =0.$ Thus,
\begin{equation*}
\Phi \left( t+a\right) =\sum_{k=1}^{\infty }A_{k}\left( t+a\right)
^{k},A_{k}=\frac{\Phi ^{\left( k\right) }\left( -a\right) }{k!}.
\end{equation*}
$\blacktriangle $

\section{The coefficients $A_{k}$}

In this section we describe a procedure that allows us to derive formulas
for computing the coefficients $A_{k},k=1,\ldots ,$ given by (\ref{at -2}).
Simultaneously we give an alternative proof of Theorem 3 which is\ based on
Theorem 1. Returning to the relation (\ref{SBorel}) rewritten as \ \
\begin{equation}
F\left( t\right) -F_{0}\left( t,r\right) =\frac{T}{2\pi i}\int_{r}^{+\infty
}e^{-\left( t+a\right) \zeta }P\left( \zeta \right) \frac{d\zeta }{\zeta },
\label{T-Laplace}
\end{equation}%
where $F_{0}\left( t,r\right) $ is given by (\ref{F0}) we represent the
integral in (\ref{T-Laplace}) in the form
\begin{equation}
\frac{T}{2\pi i}\int_{r}^{+\infty }e^{-\left( t+a\right) \zeta }P\left(
\zeta \right) \frac{d\zeta }{\zeta }=I_{n}\left( t,r\right) +J_{n}\left(
t,r\right)   \label{I+J}
\end{equation}%
where
\begin{equation}
I_{n}\left( t,r\right) =\frac{T}{2\pi i}\int_{r}^{+\infty }e^{-\left(
t+a\right) \zeta }\left( P\left( \zeta \right) -\sum_{k=0}^{n}p_{k}/\zeta
^{k}\right) \frac{d\zeta }{\zeta }  \label{I-n}
\end{equation}%
and
\begin{equation}
J_{n}\left( t,r\right) =\frac{T}{2\pi i}\sum_{k=0}^{n}p_{k}\int_{r}^{+\infty
}e^{-\left( t+a\right) \zeta }\frac{d\zeta }{\zeta ^{k+1}}.  \label{J-n-t}
\end{equation}

Given $n\in \mathbb{N}$, and expanding $I_{n}\left( t,r\right) $, $%
J_{n}\left( t,r\right) $, and $F_{0}\left( t,r\right) $ in series of the
form
\begin{equation*}
\log \left( t+a\right) \sum_{k}\beta _{k}\left( t+a\right)
^{k}+\sum_{k}\alpha _{k}\left( t+a\right) ^{k},
\end{equation*}%
we are only interested in the values of the coefficients $\alpha _{n}$ in
each series. Let us denote these coefficients by $\alpha _{n}\left( I\right)
$, $\alpha _{n}\left( J\right) $ and $\alpha _{n}\left( F_{0}\right) $ for $%
I_{n}\left( t,r\right) $, $J_{n}\left( t,r\right) $ and $F_{0}\left(
t,r\right) ,$ respectively. Then (\ref{T-Laplace}) shows that the
coefficient $A_{n}$ at $\left( t+a\right) ^{n}$ can be written as%
\begin{equation}
A_{n}=\alpha _{n}\left( I\right) +\alpha _{n}\left( J\right) +\alpha
_{n}\left( F_{0}\right) .  \label{A_n}
\end{equation}

Firstly, expanding the exponential $e^{-\left( t+a\right) \zeta }$ in (\ref%
{I-n}) into a Taylor series we have
\begin{equation}
I_{n}\left( t,r\right) =\frac{T}{2\pi i}\int_{r}^{+\infty }\sum_{m=0}^{n}%
\frac{\left( -1\right) ^{m}\left( t+a\right) ^{m}\zeta ^{m}}{m!}\left(
P\left( \zeta \right) -\sum_{k=0}^{n}p_{k}/\zeta ^{k}\right) \frac{d\zeta }{%
\zeta }+o\left( \left( t+a\right) ^{n}\right) ,t\rightarrow -a,
\label{I_n-n}
\end{equation}%
so that
\begin{equation}
\alpha _{n}\left( I\right) =\left( -1\right) ^{n}\frac{T}{2\pi i}\frac{1}{n!}%
\int_{r}^{+\infty }\left( P\left( \zeta \right) -\sum_{k=0}^{n}p_{k}/\zeta
^{k}\right) \zeta ^{n-1}d\zeta ,  \label{alpha-I}
\end{equation}%
and the integral is absolutely convergent since $\left( P\left( \zeta
\right) -\sum_{k=0}^{n}p_{k}/\zeta ^{k}\right) \zeta ^{n-1}=O\left( \frac{1}{%
\zeta ^{2}}\right) $ as $\zeta \rightarrow \infty .$ Moreover, using error
bound (\ref{errors-1}) of Theorem 1 yields the estimate for $n\in \mathbb{N}$%
\begin{equation}
\left\vert \alpha _{n}\left( I\right) \right\vert \leq \frac{\left\vert
T\right\vert M_{P}}{2\pi }\frac{\sqrt{n+3}}{a^{n}},M_{P}=\inf_{0<r<\infty
}M_{P}\left( r\right) .  \label{a-I-est}
\end{equation}

Secondly, we calculate $J_n\left( t,r\right) $\ \ using the exponential
integral
\begin{equation}
E_n\left( z\right) =\int_1^\infty e^{-zt}\frac{dt}{t^n},n=1,2,\ldots ,
\label{E-n-int}
\end{equation}
and its power series expansion
\begin{equation}
E_n\left( z\right) =\frac{\left( -z\right) ^{n-1}}{\left( n-1\right) !}%
\left( -\log z+\psi \left( n\right) \right) -\sum_{m=0,m\neq n-1}^\infty
\frac{\left( -z\right) ^m}{\left( m-n+1\right) m!},  \label{E-n}
\end{equation}
where
\begin{equation*}
\left| \arg z\right| <\pi ,\psi \left( 1\right) =-\gamma ,\psi \left(
n\right) =-\gamma +\sum_{m=1}^{n-1}\frac 1m,
\end{equation*}
and $\gamma $ is Euler's constant.

We have
\begin{equation}
J_n\left( t,r\right) =\frac T{2\pi i}\sum_{k=0}^n\frac{p_k}{r^k}%
E_{k+1}\left( r\left( t+a\right) \right)  \label{J-n}
\end{equation}
and substituting (\ref{E-n}) into (\ref{J-n}) yields
\begin{eqnarray}
J_n\left( t,r\right) &=&\frac T{2\pi i}\sum_{k=0}^np_k\frac{\left( -\left(
t+a\right) \right) ^k}{k!}\left( -\log r-\log \left( t+a\right) +\psi \left(
k+1\right) \right)  \notag \\
&&-\frac T{2\pi i}\sum_{k=0}^n\frac{p_k}{r^k}\sum_{m=0,m\neq k}^\infty \frac{%
\left( -r\left( t+a\right) \right) ^m}{\left( m-k\right) m!}.  \label{J-exp}
\end{eqnarray}

It follows from (\ref{J-exp}) that the coefficient $\alpha _{n}\left(
J\right) $ for $J_{n}\left( t,r\right) $ is given by
\begin{equation}
\alpha _{n}\left( J\right) =\frac{T}{2\pi i}\frac{\left( -1\right) ^{n}}{n!}%
p_{n}\left( -\log r+\psi \left( n+1\right) \right) -\frac{T}{2\pi i}%
\sum_{k=0}^{n-1}\frac{\left( -1\right) ^{n}r^{n-k}p_{k}}{n!\left( n-k\right)
}.  \label{alpha-J}
\end{equation}%
Using (\ref{est-p-n}) and setting $r=1,$ we have%
\begin{equation}
\left\vert \alpha _{n}\left( J\right) \right\vert =\frac{TM_{P}\left(
1\right) }{2\pi }\frac{\psi \left( n+1\right) }{a^{n}}\left( 1+o\left(
1\right) \right) ,n\rightarrow \infty .  \label{a-J-est}
\end{equation}

Finally, it follows from (\ref{F-0-r}) that
\begin{equation}
\alpha _{n}\left( F_{0}\right) =\frac{1}{2\pi i}\frac{1}{n!}\int_{C\left(
r\right) }e^{-a\zeta }P\left( \zeta \right) \frac{d\zeta }{\zeta },
\label{alpha-F-0}
\end{equation}%
and%
\begin{equation}
\left\vert \alpha _{n}\left( F_{0}\right) \right\vert \leq \frac{e^{ar}M_{P}%
}{n!}  \label{a-F_0-est}
\end{equation}%
So we have proved that the coefficient $A_{n}$ at $\left( t+a\right) ^{n}$
is represented in the form given by (\ref{A_n}), where $\alpha _{n}\left(
I\right) ,\alpha _{n}\left( J\right) ,$ and $\alpha _{n}\left( F_{0}\right) $
are given by (\ref{alpha-I}), (\ref{alpha-J}), and (\ref{alpha-F-0}),
respectively. Analysis of the estimates for $\alpha _{n}\left( I\right)
,\alpha _{n}\left( J\right) ,$ and $\alpha _{n}\left( F_{0}\right) $ given
by (\ref{a-I-est}), (\ref{a-J-est}), and (\ref{a-F_0-est}) shows that the
term (\ref{alpha-I}) gives the main contribution to the asymptotics of $A_{n}
$ as $n\rightarrow \infty .$ We have the estimae%
\begin{equation}
\left\vert A_{n}\right\vert =\frac{TM_{P}\left( 1\right) }{2\pi }\frac{\sqrt{%
n+3}}{a^{n}}\left( 1+o\left( 1\right) \right) ,n\rightarrow \infty
\label{A-est}
\end{equation}%
which shows that the power series $\sum_{n=0}^{\infty }A_{n}t^{n}$ is
absolutely convergent in the circle of radius $a.$

It turns out that we have an alternative proof of Theorem 3 based on Theorem
1. Indeed, combining expressions (\ref{T-Laplace}), (\ref{I+J}), (\ref{I_n-n}%
), (\ref{J-exp}), we have for $\left\vert \arg \left( t+a\right) \right\vert
<a$%
\begin{equation}
F\left( t\right) =-\frac{T}{2\pi i}\sum_{k=0}^{n}\frac{p_{k}}{r^{k}}\frac{%
\left( -r\left( t+a\right) \right) ^{k}}{k!}\log \left( t+a\right)
+\sum_{n=0}^{n}A_{n}\left( t+a\right) ^{n}+o\left( \left( t+a\right)
^{n}\right) ,t\rightarrow -a,  \label{F-appr-a}
\end{equation}%
where $A_{n}$ are given by (\ref{A_n}). Since $p_{k}=F^{\left( k\right)
}\left( 0\right) $ we have
\begin{equation*}
-\lim_{n\rightarrow \infty }\frac{T}{2\pi i}\sum_{k=0}^{n}\frac{p_{k}}{r^{k}}%
\frac{\left( -r\left( t+a\right) \right) ^{k}}{k!}\log \left( t+a\right) =-%
\frac{T}{2\pi i}F\left( -\left( t+a\right) \right) \log \left( t+a\right) ,
\end{equation*}%
so (\ref{F-appr-a}) can be rewritten as%
\begin{equation*}
F\left( t\right) =-\frac{T}{2\pi i}F\left( -\left( t+a\right) \right) \log
\left( t+a\right) +\sum_{n=0}^{n}A_{n}\left( t+a\right) ^{n}+o\left( \left(
t+a\right) ^{n}\right) ,t\rightarrow -a.
\end{equation*}

It remains only to evaluate more accurately the remainder $o\left( \left(
t+a\right) ^{n}\right) $ in (\ref{I_n-n}). Let us return to the expression (%
\ref{I-n}) and represent it in the form%
\begin{equation*}
I_{n}\left( t,r\right) =\frac{T}{2\pi i}\int_{r}^{+\infty }e^{-\left(
t+a\right) \zeta }\left( P\left( \zeta \right) -\sum_{k=0}^{n}p_{k}/\zeta
^{k}\right) \frac{d\zeta }{\zeta }
\end{equation*}%
which can be rewritten as%
\begin{equation}
I_{n}\left( t,r\right) =\frac{T}{2\pi i}\int_{r}^{+\infty }e^{-\left(
t+a\right) \zeta }\mathrm{P}_{n+1}\left( \zeta \right) \frac{d\zeta }{\zeta }%
,  \label{I-last}
\end{equation}%
where $\mathrm{P}_{n+1}\left( \zeta \right) $ is given by (\ref{remainders}%
). Using the expansion%
\begin{equation*}
e^{-\left( t+a\right) \zeta }=\sum_{m=0}^{n-1}\frac{\left( -1\right)
^{m}\left( t+a\right) ^{m}\zeta ^{m}}{m!}+\frac{\left( -1\right) ^{n}\left(
t+a\right) ^{n}\zeta ^{n}}{n!}e^{-\left( t+a\right) \zeta \left( r\right) },
\end{equation*}%
where $r<\zeta \left( r\right) <\infty ,$ and notation (\ref{alpha-I}), the
expression (\ref{I-last}) can be rewritten as%
\begin{equation*}
I_{n}\left( t,r\right) =\sum_{k=0}^{n-1}\alpha _{k}\left( I\right) \left(
t+a\right) ^{k}+R\left( t,n,r\right) ,
\end{equation*}%
where%
\begin{equation}
R\left( t,n,r\right) =\frac{T}{2\pi i}\frac{\left( -1\right) ^{n}\left(
t+a\right) ^{n}}{n!}\int_{r}^{+\infty }\zeta ^{n}e^{-\left( t+a\right) \zeta
\left( r\right) }\mathrm{P}_{n+1}\left( \zeta \right) \frac{d\zeta }{\zeta }.
\label{R-t-n}
\end{equation}%
Applying the estimate for the remainders given by (\ref{errors-1}) we have
for $n=0,1,\ldots ,$%
\begin{equation*}
\left\vert \mathrm{P}_{n+1}\left( \zeta \right) \right\vert \leq M_{P}\left(
r\right) \frac{\left( n+1\right) !\sqrt{n+4}}{a^{n+1}\zeta ^{n+1}}{,}
\end{equation*}
which together with (\ref{R-t-n}) yields the inequality%
\begin{equation*}
\left\vert R\left( t,n,r\right) \right\vert \leq \frac{TM_{P}\left( r\right)
}{2\pi }\frac{\left( n+1\right) \sqrt{n+4}}{a^{n}}\left\vert t+a\right\vert
^{n},n=0,1,\ldots .
\end{equation*}%
Using estimate (\ref{A_n}) this completes the alternative proof of Theorem
3. $\blacktriangle $

\section{Appendix. The perturbed Whittaker equation}

We have replaced the differential equation (\ref{T-Airy}) by a system of
functional equations
\begin{eqnarray*}
P_1\left( \zeta e^{\pi i}\right) &=&P_1\left( \zeta e^{-\pi i}\right)
+Te^{-a\zeta }P_1\left( \zeta \right) , \\
P_2\left( \zeta e^{\pi i}\right) &=&P_2\left( \zeta e^{-\pi i}\right)
+Te^{a\zeta }P_2\left( \zeta \right) .
\end{eqnarray*}

We note that these two equations of the system are not linked. A reason for
this is that the coefficient $A\left( \zeta \right) $ is an even function.
However, typically, such equations are intertwined. It is enough to add the
term $\frac{\mathfrak{b}}\zeta $ to $A\left( \zeta \right) $ to obtain from (%
\ref{T-Airy}) an intertwined differential equation
\begin{equation*}
\frac{d^2u}{d\zeta ^2}=\left( \frac{a^2}4+\frac b\zeta +\frac{a_0}{\zeta ^2}+%
\frac{a_1}{\zeta ^4}+\ldots \right) u.
\end{equation*}
We consider a more general equation
\begin{equation}
\frac{d^2u}{d\zeta ^2}=\left( \frac{a^2}4+\frac b\zeta +\frac{B\left( \zeta
\right) }{\zeta ^2}\right) u,  \label{p-Wh}
\end{equation}
where $B\left( \zeta \right) =\sum_{k=0}^\infty b_k/\zeta ^k$ is an entire
function of $1/\zeta $. A special case $a=1,b=-\kappa ,b_0=\mu ^2-1/4,b_k=0$
for $k=1,2,\ldots ,$ is known as Whittaker's differential equation. Thus, (%
\ref{p-Wh}) can be considered as a perturbed Whittaker equation. Let us
assume that $a>0$. Then it can be proved that there exists a pair of
linearly independent solutions $u_1\left( \zeta \right) $ and $u_2\left(
\zeta \right) $ of (\ref{p-Wh}) such that
\begin{equation}
\begin{array}{l}
u_1\left( \zeta \right) =e^{-\frac a2\zeta }\zeta ^{-\frac ba}P_1\left(
\zeta \right) ,-\frac{3\pi }2<\arg \zeta <\frac{3\pi }2, \\
u_2\left( \zeta \right) =e^{\frac a2\zeta }\zeta ^{\frac ba}P_2\left( \zeta
\right) ,-\frac \pi 2<\arg \zeta <\frac{5\pi }2%
\end{array}
\label{2}
\end{equation}
where
\begin{equation}
P_1\left( \zeta \right) ,P_2\left( \zeta \right) =1+o\left( 1\right)
\label{2-asymp}
\end{equation}
as $\zeta \rightarrow \infty $ along any ray of these sectorial regions. It
can be shown that the solutions $u_1\left( \zeta \right) $ and $u_2\left(
\zeta \right) $ are uniquely determined by their asymptotics given by (\ref%
{2}) and (\ref{2-asymp}) and that $P_1\left( \zeta \right) $ and $P_2\left(
\zeta \right) $ admit analytical continuation along any path not crossing $%
\zeta =0$.

Since $u_1\left( \zeta e^{2\pi i}\right) $ and $u_2\left( \zeta e^{2\pi
i}\right) $ are also solutions of (\ref{p-Wh}) and since $u_1\left( \zeta
\right) $ and $u_2\left( \zeta e^{2\pi i}\right) $ are linearly independent
solutions we have
\begin{equation}
\begin{array}{l}
u_1\left( \zeta e^{2\pi i}\right) =Au_1\left( \zeta \right) +Bu_2\left(
\zeta e^{2\pi i}\right) \\
u_2\left( \zeta e^{2\pi i}\right) =Cu_2\left( \zeta \right) +Du_1\left(
\zeta \right)%
\end{array}
,  \label{3}
\end{equation}
where $A,B,C,D$ are complex constants. Using (\ref{2}), we have
\begin{equation*}
\begin{array}{c}
u_1\left( \zeta e^{2\pi i}\right) =e^{-2\pi i\frac ba}e^{-\frac a2\zeta
}\zeta ^{-\frac ba}P_1\left( \zeta e^{2\pi i}\right) \\
u_2\left( \zeta e^{2\pi i}\right) =e^{2\pi i\frac ba}e^{\frac a2\zeta }\zeta
^{\frac ba}P_2\left( \zeta e^{2\pi i}\right)%
\end{array}
,
\end{equation*}
which allows us to rewrite (\ref{3}) in the form
\begin{equation}
\begin{array}{l}
e^{-2\pi i\frac ba}e^{-\frac a2\zeta }\zeta ^{-\frac ba}P_1\left( \zeta
e^{2\pi i}\right) =Ae^{-\frac a2\zeta }\zeta ^{-\frac ba}P_1\left( \zeta
\right) +Be^{2\pi i\frac ba}e^{\frac a2\zeta }\zeta ^{\frac ba}P_2\left(
\zeta e^{2\pi i}\right) \\
e^{2\pi i\frac ba}e^{\frac a2\zeta }\zeta ^{\frac ba}P_2\left( \zeta e^{2\pi
i}\right) =Ce^{\frac a2\zeta }\zeta ^{\frac ba}P_2\left( \zeta \right)
+De^{-\frac a2\zeta }\zeta ^{-\frac ba}P_1\left( \zeta \right)%
\end{array}
.  \label{4}
\end{equation}

Let us simplify (\ref{4})
\begin{eqnarray}
e^{-2\pi i\frac ba}P_1\left( \zeta e^{2\pi i}\right) &=&AP_1\left( \zeta
\right) +Be^{2\pi i\frac ba}e^{a\zeta }\zeta ^{\frac{2b}a}P_2\left( \zeta
e^{2\pi i}\right) ,  \label{5.1} \\
e^{2\pi i\frac ba}P_2\left( \zeta e^{2\pi i}\right) &=&CP_2\left( \zeta
\right) +De^{-a\zeta }\zeta ^{-\frac{2b}a}P_1\left( \zeta \right) .
\label{5.2}
\end{eqnarray}
and analyze the last pair of equations. Considering the equation (\ref{5.1})
we assume that $\arg \zeta \in \left( -\frac{3\pi }2,-\frac \pi 2\right) ,$
thus $\arg \left( \zeta e^{2\pi i}\right) \in \left( \frac \pi 2,\frac{3\pi }%
2\right) .$ Then letting $\zeta $ go to infinity, and taking into account
that $P_1\left( \zeta \right) \rightarrow 1,P_2\left( \zeta e^{2\pi
i}\right) \rightarrow 1,$and that the term $Be^{a\zeta }\zeta ^{\frac{2b}%
a}P_2\left( \zeta e^{2\pi i}\right) $ is exponentially small, it follows
that
\begin{equation}
A=e^{-2\pi i\frac ba}.  \label{6}
\end{equation}
A similar analysis for (\ref{5.2}) and for $\arg \zeta \in \left( -\frac \pi
2,\frac \pi 2\right) $\ shows that
\begin{equation}
C=e^{2\pi i\frac ba}.  \label{7}
\end{equation}
Using (\ref{6}) and (\ref{7}), the system of equations (\ref{5.1}) and (\ref%
{5.2}) can be rewritten as
\begin{eqnarray}
P_1\left( \zeta e^{2\pi i}\right) &=&P_1\left( \zeta \right) +Be^{4\pi
i\frac ba}e^{a\zeta }\zeta ^{\frac{2b}a}P_2\left( \zeta e^{2\pi i}\right) ,
\label{8.1} \\
P_2\left( \zeta e^{2\pi i}\right) &=&P_2\left( \zeta \right) +De^{-2\pi
i\frac ba}e^{-a\zeta }\zeta ^{-\frac{2b}a}P_1\left( \zeta \right) .
\label{8.2}
\end{eqnarray}
Setting
\begin{equation*}
T_1=Be^{4\pi i\frac ba},T_2=De^{-2\pi i\frac ba}
\end{equation*}
we have finally
\begin{eqnarray}
P_1\left( \zeta e^{2\pi i}\right) &=&P_1\left( \zeta \right) +T_1e^{a\zeta
}\zeta ^{\frac{2b}a}P_2\left( \zeta e^{2\pi i}\right) ,  \label{P-1} \\
P_2\left( \zeta e^{2\pi i}\right) &=&P_2\left( \zeta \right) +T_2e^{-a\zeta
}\zeta ^{-\frac{2b}a}P_1\left( \zeta \right) .  \label{P-2}
\end{eqnarray}
Note that the above pair of relations can also be rewritten as
\begin{eqnarray}
P_1\left( \zeta e^{\pi i}\right) &=&P_1\left( \zeta e^{-\pi i}\right)
+T_1e^{-2\pi i\frac ba}e^{-a\zeta }\zeta ^{\frac{2b}a}P_2\left( \zeta e^{\pi
i}\right) ,  \label{9.1} \\
P_2\left( \zeta e^{\pi i}\right) &=&P_2\left( \zeta e^{-\pi i}\right)
+T_2e^{2\pi i\frac ba}e^{a\zeta }\zeta ^{-\frac{2b}a}P_1\left( \zeta e^{-\pi
i}\right) .  \label{9.2}
\end{eqnarray}

Our principal idea is to consider the system of equations (\ref{9.1}) and (%
\ref{9.2}) separately of differential equation. Assuming that $P_1\left(
\zeta \right) $ and $P_2\left( \zeta \right) $ are analytic and bounded in
the regions $-\frac{3\pi }2<\arg \zeta <\frac{3\pi }2$ and $-\frac \pi
2<\arg \zeta <\frac{5\pi }2,$ respectively and representing these functions
in the form
\begin{equation*}
P_1\left( \zeta \right) =\zeta \int_0^\infty e^{-\zeta t}F_1\left( t\right)
dt
\end{equation*}
and
\begin{equation*}
P_2\left( \zeta \right) =\zeta \int_0^{\infty \cdot e^{\pi i}}e^{-\zeta
t}F_2\left( t\right) dt,
\end{equation*}
we claim the following statement.

\textbf{Theorem. }\textit{(i) Functions }$F_{1}\left( t\right) $\textit{\
and }$F_{2}\left( t\right) $\textit{\ admit analytical continuation to the }$%
t$\textit{-plane punctured at points }$t=0$\textit{\ and }$t=-a$\textit{\
for }$F_{1}\left( t\right) $\textit{\ and at points }$t=0$\textit{\ and }$%
t=a $\textit{\ for }$F_{2}\left( t\right) ;$

\textit{(ii) there exist branches of }$F_{1}\left( t\right) $\textit{\ and }$%
F_{2}\left( t\right) $\textit{\ analytic in the }$t$\textit{-plane cut along
}$\left( -\infty ,-a\right) $\textit{\ and }$\left( a,+\infty \right) , $%
\textit{\ respectively;}

\textit{(iii) }$F_{1}\left( t\right) $\textit{\ and }$F_{2}\left( t\right) $%
\textit{\ satisfy a dual system of monodromic relations.}

The most nontrivial is the assertion (ii). Setting
\begin{equation*}
p_k^{\left( 1\right) }=F_1^{\left( k\right) }\left( 0\right) ,p_k^{\left(
2\right) }=F_2^{\left( k\right) }\left( 0\right) ,k=0,1\ldots ,
\end{equation*}
it follows from this assertion, in particular, that
\begin{align*}
-\frac{3\pi }2& <\arg \zeta <\frac{3\pi }2\Rightarrow \lim_{\zeta
\rightarrow \infty }P_1\left( \zeta \right) =F_1\left( 0\right) , \\
-\frac \pi 2& <\arg \zeta <\frac{5\pi }2\Rightarrow \lim_{\zeta \rightarrow
\infty }P_2\left( \zeta \right) =F_2\left( 0\right) ,
\end{align*}
and moreover that $P_1\left( \zeta \right) $ and $P_2\left( \zeta \right) $
can be expanded into asymptotic series
\begin{align*}
P_1\left( \zeta \right) & \sim \sum_{k=0}^\infty \frac{p_k^{\left( 1\right) }%
}{\zeta ^k},-\frac{3\pi }2<\arg \zeta <\frac{3\pi }2, \\
P_2\left( \zeta \right) & \sim \sum_{k=0}^\infty \frac{p_k^{\left( 2\right) }%
}{\zeta ^k},-\frac \pi 2<\arg \zeta <\frac{5\pi }2.
\end{align*}

\section{Conclusion}

Our aim for the future is to extend the approach described here to more
general systems of functional monodromic equations that are generated by
linear differential equations or systems of differential equations with an
irregular singular point of arbitrary Poincar\'{e} rank at infinity. The
pBde is the simplest differential equation that can be reduced to a single
functional equation while preserving most of the difficulties arising in the
general case. It is for this reason that in this initial study we have
limited our attention to a detailed consideration of the pBde.

\section{Acknowledgment}

Authors wish to express their gratitude to Sergey Suslov for his useful
suggestions.

{Mathematics 38, \textit{FEIS},\newline
Swinburne University of Technology\newline
PO Box 218 Hawthorn 3122;\newline
e-mail: vgurarii@swin.edu.au,}

School of Information Technology and Mathematical Sciences\newline
University of Ballarat\newline
P.O. Box 663 Ballarat, 3353, \newline
e-mail: dgillam@ballarat.edu.au \newline
VIC\newline
Australia

VIC\newline
Australia.

\end{document}